\documentclass[a4paper, 12pt]{amsproc}
\usepackage{mathrsfs}
\usepackage{amssymb}
\usepackage{amsmath}
\title[Infinite Dimensional Representations of Reductive Groups]
{Some Infinite Dimensional Representations of Reductive Groups with Frobenius Maps}

\author[Nanhua Xi]{ Nanhua XI$^{*}$}
\address{$^{*}$
Institute of Mathematics\\
Chinese Academy of Sciences\\
Beijing, 100190\\
China, and School of Mathematical Sciences, University of Chinese
Academy of Sciences, Beijing, China} \email{nanhua@math.ac.cn}

\thanks{Supported in part by a grant of the National Natural Science Foundation of China (No. 11321101).}

\begin{document}
\baselineskip=18pt
\begin{abstract}
In this paper we  construct certain irreducible infinite dimensional representations of algebraic groups
with Frobenius maps. In particular, a few classical results of Steinberg and Deligne \& Lusztig on
complex representations of finite groups of Lie type are extended to reductive algebraic groups with Frobenius
maps.
\end{abstract}

\maketitle

\def\Cal{\mathcal}
\def\bold{\mathbf}
\def\ca{\mathcal A}
\def\cdz{\mathcal D_0}
\def\cd{\mathcal D}
\def\cdo{\mathcal D_1}
\def\bold{\mathbf}
\def\l{\lambda}
\def\le{\leq}

The construction of induced representations of Frobenius for finite groups has various generalizations for
infinite groups. It seems that for infinite groups of Lie type the original form of construction of Frobenius
was not used much. In this paper, we try to study abstract representations of algebraic groups by using the
original construction of Frobenius directly. We are mainly interested in reductive groups with Frobenius
maps. A few classical results of Steinberg and Deligne \& Lusztig on complex representations of finite
groups of Lie type are extended to reductive algebraic groups with Frobenius maps, see Propositions 2.3
and 2.4, Theorems 3.2 and 3.4, etc..

The paper is organized as follows. In Section 1 we give some trivial extensions for several results in representation
theory of finite groups and introduce the concept of quasi-finite groups (see Subsection 1.8).
A few general results on irreducibility of a representation for a quasi-finite group are established, if the
representation is a ¡°limit¡± of the irreducible representations (Lemmas 1.5 and 1.6). A partial generalization
of Mackeys criterion on irreducibility of induced modules for quasi-finite groups is given (see
Subsection 1.9).

In Section 2 we consider algebraic groups with split BN-pairs. The main objects of this section are
induced representations of certain one-dimensional representations of a Borel subgroup of an algebraic
group with split BN-pairs. In particular, the Steinberg module of the algebraic group is constructed
(Proposition 2.3(b)). In Section 3 we consider reductive groups with Frobenius maps. The main results
are Theorems 3.2 and 3.4. The first one says that the Steinberg module of a reductive group over an
algebraically closed field of positive characteristic is irreducible when the base field of the Steinberg
module is the field of complex numbers or the ground field of the reductive group, the second one says
that the induced representations of certain one-dimensional complex representations of a Borel subgroup
are irreducible.

Gelfand-Graev modules of reductive groups are defined in Section 4, which are similar to those for finite
groups of Lie type. In Section 5 a few questions are raised. In Section 6 we discuss type $A_1$. Section 7 is
devoted to discussing representations of some infinite Coxeter groups and infinite dimensional groups of
Lie type.

This work was partially motivated by trying to find an algebraic counterpart for Lusztig¡¯s theory of
character sheaves, the author is grateful to Professor G. Lusztig for his series of lectures on character
sheaves delivered at the Academy of Mathematics and Systems Science, Chinese Academy of Sciences,
Beijing, in 2012.

\section{General Setting}
\def\ind{{\text {Ind}}}
\def\Hom{{\text {Hom}}}
\def\res{{\text {Res}}}

\subsection{} In this section we give some trivial
extensions for several results in representation theory of
finite groups. There are many good references, say  [S] and [CR].

Let $H$ be a subgroup of a group $G$ and $k$ a field. In this section all modules
are assumed to be over $k$. For an $H$-module $M$, we can consider the naive  induced module of $M$:
$$\ind_H^GM=k G\otimes_{k H}M,$$
where $k G$ and $k H$ are the group algebras of $G$ and $H$ over
the field $k$ respectively.

As the case of finite groups, we can define the induced module in
another way.   Let $\mathscr{M}$ be the set of all functions
$f:G\to M$ satisfying $f(gh)=h^{-1}f(g)$ for any $g\in G$ and
$h\in H$.  For $f\in\mathscr{M}$ and $x\in G$, set
$xf(g)=f(x^{-1}g)$. This defines a $kG$-module structure on
$\mathscr M$.

Let $G/H$ be the set of left cosets of $H$ in $G$. A function $f:G
\to M$ is called to have finite support on $G/H$ if all $f(gH)=0$
except for finitely many left cosets of $H$ in $G$. Let
$\mathscr{M}_0$ be the subset of $\mathscr{M}$ consisting of all
functions $f$ in $\mathscr{M}$ with finite support on $G/H$. It is
clear that $\mathscr{M}_0$ is a $kG$-submodule of $\mathscr{M}.$
The following result is known for finite groups.

\subsection{Lemma.} The $kG$-module $\mathscr{M}_0$ is isomorphic to
the induced module $kG\otimes_{kH}M$.

Proof. Let $\{g_i\}_{i\in I}$ be a set of representatives of left
cosets of $H$ in $G$. Then the map $f\to \sum_i g_i\otimes f(g_i)$
defines an isomorphism of $kG$-module from $\mathscr{M}_0$ to
$kG\otimes_{kH}M$.

\medskip

The induced modules above are extremely important in representation
theory of finite groups and Lie algebras, but seem not studied much
for infinite groups of Lie type. We have some trivial properties for
these induced modules, such as Frobenius's reciprocity, etc..

\subsection{Lemma.} (a) Let $M$ be an $H$-module and $N$ be a
$G$-module. Then we have
 $$\Hom_G(\ind_H^GM,N)\simeq \Hom_H(M,\res_HN)$$
 and
 $$\ind_H^G(M\otimes \res_HN)\simeq \ind_H^GM\otimes N.$$
 where $\res_H$ denotes the restriction functor from $G$-modules to $H$-modules.

(b) Let $H\subset K$ be subgroups of $G$ and $M$ an $H$-module, then
$\ind_K^G(\ind_H^KM)$ is isomorphic to $\ind_H^GM$.

The following result should be known.

\subsection{Lemma.} Assume that $G$ is commutative and each element of $G$ has finite order. If $k$ is
algebraically closed, then any irreducible representation of $G$
over $k$ is one dimensional.

Proof. Let $M$ be an irreducible $k G$-module.  By Schur lemma,
End${}_GM$ is a division algebra over $k$. Since $G$ is commutative,
for any $x\in G$, the map $\varphi_x:\ M\to M,\ a\to xa$ is a
$G$-homomorphism. Since $x$ has finite order, $\varphi_x$ is algebraic
over $k$. Now $k$ is algebraically closed, so $\varphi_x$ must be in
$k$, that is, $x$ acts on $M$ by multiplication of a scalar in $k$.
 Thus $M$ must be one-dimensional since $M$ is
irreducible. The lemma is proved.

{\bf Remark:} I am grateful to Binyong Sun for pointing out that the lemma above can not be extended to arbitrary commutative groups, say, the field $\mathbb{C}$ is an irreducible representation of $\mathbb{C}^*$, but it is of infinite dimension over $\bar {\mathbb{Q}}$.

\def\bbf{\mathbb{F}_q}
 \def\bbfa{\mathbb{F}_{q^a}}
 \def\bbbf{\bar{\mathbb{F}}_q}
 \def\bbc{\mathbb{C}}

 \subsection{Lemma.} Let $(I, \preceq)$ be a directed set, $\{A_i,\ f_{ij}\}$ be a direct system of algebras
 over $k$ (resp. groups) and $\{M_i,\ \varphi_{ij}\}$ be a direct system of vector spaces over $k$.
 Assume that $M_i$ is $A_i$-module for each $i$ and for any $i,j\in
 I$ with $i\preceq j$, the  homomorphism $\varphi_{ij}:M_i\to M_j$ is compatible
 with the homomorphism $f_{ij}: A_i\to A_j$, that is, $\varphi_{ij}(ux)=f_{ij}(u)\varphi_{ij}(x)$
 for any $u\in A_i$ and $x\in M_i$. Then $M=\underset{\rightarrow}{\lim} M_i$ is naturally a module of $A=\underset{\rightarrow}{\lim} A_i$. Moreover, if $M_i$ is
 irreducible $A_i$-module for each $i\in I$, then $M$ is irreducible
 $A$-module.

 Proof. The proof is easy. For convenience, we give the details. Let $\sqcup_i M_i$ be the disjoint union of all $M_i$.
 By definition, $M=\sqcup_i M_i /\thicksim$, here for $x'\in M_i$
 and $y'\in M_j$, $x'\thicksim y'$ if and only if there exists some $r\in I$ such that $\varphi_{ir}(x')=\varphi_{jr}(y')$.

 Let $x\in M$ (resp. $u\in A$). Choose $x'\in M_i$  (resp. $u'\in A_j$) such that $x'$ (resp. $u'$) is in the equivalence class $x$ (resp. $u$). Choose
 $r\in I$ such that $i\preceq r$ and $j\preceq r$. Then set $ux$ to be the
 class containing $f_{jr}(u')\varphi_{ir}(x').$  One can verify that this defines an $A$-module structure on $M$.

 Assume that $M_i$ is irreducible $A_i$-module for each $i$. To prove that $M$ is irreducible $A$-module it suffices to
  prove that  $M=Ax$ for any nonzero element $x$ in
 $M$.  Assume that $ x$ and $ y$ are two nonzero elements in $M$. Let $x'\in M_i$
(resp. $y'\in M_j$) be an element in the equivalence class $x$
(resp. $y$). Choose $r\in I$ such that $i\preceq r$ and $j\preceq
r$. Then $x''=\varphi_{ir}(x')$ (resp. $y''=\varphi_{jr}(y')$) is in
the class $x$ (resp. $y$). Since $M_r$ is irreducible $A_r$-module,
there exists $u'\in A_r$ such that $x''=u'y''$. Let $u$ be the
equivalence class containing $u'$. Then $u$ is element of $A$ and
$ux=y$. The lemma is proved.

\medskip

 The following two simple lemmas  will be used frequently.

\medskip

\subsection{Lemma.} (a) Let $A$ be an algebra over $k$ and $M$ be
an
 $A$-module. Assume that $A$ has a sequence of
 subalgebras $A_1,\ A_2,\ ,...,\ A_n,\ ... $
 and   that $M$ has a sequence of $k$-subspaces
 $M_1,\ M_2,\ ...,\ M_n,\ ...$ such that $M$ is the union of all $M_i$ and for any positive integers
 $i,j$ there exists positive integer $r$
 such that $M_i$ and $M_j$ are contained in $M_r$. If  $M_i$ is an irreducible $A_i$-submodule of $M$ for any $i$,
 then $M$ is an irreducible
 $A$-module.

 (b) Let $G$ be a group and $M$ be a $G$-module. Assume that $G$ has a sequence of subgroups
 $G_1,\ G_2,\ ,...,\ G_n,\ ... $
 and   that $M$ has a sequence of $k$-subspaces
 $M_1,\ M_2,\ ...,\ M_n,\ ...$ such that $M$ is the union of all $M_i$ and for any positive integers $i,j$ there exists integer $r$
 such that $M_i$ and $M_j$ are contained in $M_r$. If  $M_i$ is an irreducible $G_i$-submodule of $M$ for any $i$,
 then $M$ is an irreducible
 $G$-module.

 Proof. (a) We only need to prove that $M=Ax$ for any nonzero element $x$ in
 $M$. Assume that $x$ and $y$ are two nonzero elements in $M$. Then
 we can find some positive integer $i$ such that both $x$ and $y$ are contained in
 $M_i$. Since $M_i$ is irreducible $A_i$-submodule of $M$, there exists $u$ in
 $A_i$ such that $ux=y$. Therefore $M=Ax$.

 (b) Applying (a) to $A=k G$ and $A_i=k[G_i]$, we see that (b) is a
 special case of (a).

 The lemma is proved.

\subsection{Lemma.} Let $H$ be a subgroup of $G$ and $M$ a $k H$-module.
  Assume that $G$ has a sequence  $G_1,\ G_2,\ ,...,\ G_n,\ ...$ of subgroups
 such that $G$ is the union of all $G_i$ and for any positive integers $i,j$ there exists integer $r$
 such that $G_i$ and $G_j$ are contained in $G_r$. Then the following
 results hold.

(a) As $G_i$-modules, $kG_i\otimes_{k(G_i\cap H)}M$ is isomorphic
to the $G_i$-submodule $Y_i$ of $k G\otimes_{k H}M$ generated by
all $x\otimes m$, where $x\in kG_i$ and $m\in M$.

(b) $k G\otimes_{kH}M$ is the union of all $Y_i$.

(c) $k G\otimes_{k H}M$ is irreducible if each $Y_i$ is
irreducible $G_i$-module.

(d) Let  $M_i$ be an $H_i=H\cap G_i$-submodule of $M$. Then we have
a natural homomorphism of $G_i$-module $\varphi_i:
kG_i\otimes_{kH_i}M_i\to k G\otimes_{k H}M$. If $M$ is the union of
all $M_i$ and $M_i$ is a subspace of $M_j$ whenever $G_i$ is a
subgroup of $G_j$, then $k G\otimes_{k H}M$ is the union of all the
images Im$\varphi_i$.

Proof. (a), (b) and (d) are clear, (c) follows (b) and Lemma 1.6
(b).

\subsection{}   Let $A$ be an algebra over a field $k$. Assume
that $A$ has a sequence of subalgebras $A_1,\ A_2,\ ,...,\ A_n,\
..., $
 such that $A$ is the union of all $A_i$ and for any positive integers $i,j$ there exists integer $r$
 such that $A_i$ and $A_j$ are contained in $A_r$.  We can consider a category $\mathscr{F}$ of
$A$-modules whose objects are those $A$-modules $M$ with a finite
dimensional $A_i$-submodule $M_i$ for each $i$ such that $M$ is the
union of all $M_i$ and for any positive integers $i,j$,  $M_i$ and
$M_j$ are contained in $M_r$ whenever both $A_i$ and $A_j$ are
contained in $A_r$. Let $(M,M_i)$ and $(N,N_i)$ be two objects in  $\mathscr{F}$. The morphisms from $(M,M_i)$ to $(N,N_i)$ are just those homomorphisms of $A$-module from $M$ to $N$ such that $f(M_i)\subset N_i$ for all $i$. Clearly $\mathscr{F}$ is an
abelian category.

Let $A$  be as above.  We say that  $A$   is quasi-finite if all
$A_i$ are finite dimensional over $k$. Similarly we say that a group
$G$ is quasi-finite if  $G$ has a sequence $G_1,\ G_2,\ ,...,\ G_n,\
... $ of finite subgroups
 such that $G$ is the union of all $G_i$ and for any positive integers $i,j$ there exists integer $r$
 such that $G_i$ and $G_j$ are contained in $G_r$.  The sequence $G_1,G_2,G_3,...$ is called a quasi-finite sequence of $G$.
 A
subgroup of a quasi-finite group is clearly quasi-finite. Clearly if
a group $G$ is quasi-finite then the group algebra $k G$ is a
quasi-finite algebra over $k$.

{\bf Example.} (1) Let $W_n$ be a Weyl group of one type $A_n$
(resp. $\ B_n\ (n\ge 2),\ D_n\ (n\ge 4))$, then we have a canonical
imbedding $W_n\to W_{n+1}$. Let $W_{\infty}=\cup_{n} W_n$. Then
$W_{\infty}$ is a quasi-finite group and is also a Coxeter group.

(2) Let $\bbf$ be a finite field of $q$ elements and $\bbbf$ its
algebraic closure.  The additive group of $\bbbf$ is quasi-finite
and is the union of all $\mathbb{F}_{q^{a}}$, $a=1,2,...$. Also the
multiplication group $\bbbf^*$ is quasi-finite and is the union of
all $\mathbb{F}^*_{q^{a}}$, $a=1,2,\ ...$.

(3) Let $G$ be an algebraic group defined over $\bbf$. By (2) we see
that the $\bbbf$-points $G(\bbbf)$ of $G$ is quasi-finite and is the
union of all $G(q^{a}),\ a\ge 1$, where $G(q^a)$ is the
$\mathbb{F}_{q^a}$-points of $G$.

(4) Let $G_n$ be  $GL_n(k)$ (resp.  $SL_n(k)$, $SO_{2n}(k),$
$SO_{2n+1}(k),$ $Sp_{2n}(k)$). Then $G_n$ is naturally embedded
into $G_{n+1}$. Let $G_{\infty}$ be the union of all $G_n$. If $k$
is finite then $G_\infty$ is quasi-finite.

More generally, direct union of quasi-finite groups is also
quasi-finite, in particular, $G_\infty$ is quasi-finite if
$k=\bar{\mathbb{F}}_q$. (I am grateful to a referee for pointing
out this  fact.)

{\it In the rest of this section we assume that all groups are
quasi-finite unless other specifications are given. For a
quasi-finite group $G$, we fix a quasi-finite sequence
$G_1,G_2,G_3,...$. For a subgroup $H$ of $G$, the quasi-finite
sequence of $H$ is chosen to be $H\cap G_1,H\cap G_2,H\cap G_3,$...,
called the quasi-finite sequence of $H$ induced from the given
quasi-finite sequence of $G$.}

Assume that $N$ is a finitely generated $G$-module, say generated by
$x_1,...,x_n$. For each positive integer $i$, let $N_i$ be the
$G_i$-submodule of $M$ generated by $x_1,...,x_n$. Then $N_i$ is a
subspace of $N_j$ if $G_i$ is a subgroup of $G_j$, and  $N$ is the
union of all $N_i$.

We shall say that an irreducible module (or representation) $N$ of
$G$ is {\it quasi-finite} (with respect to the quasi-finite sequence
$G_1,G_2,G_3,...)$ if it has a sequence of subspaces $N_1$, $ N_2$,
$ N_3, $ ... of $N$ such that (1) each $N_i$ is an irreducible
$G_i$-submodule of $N$, (2) if $G_i$ is a subgroup of $G_j$, then
$N_i$ is a subspace of $N_j$, (3) $N$ is the union of all $N_i$. The
sequence $N_1$, $ N_2$, $ N_3, $ ... will be called a quasi-finite
sequence of $N$. If the intersection $\cap_i N_i$ of all $N_i$ is
nonzero, then a nonzero element in the intersection $\cap_i N_i$
will be called primitive since such an element generates an
irreducible $G_i$-submodule of $N$ for any $i$. It is often that
$G_1$ is a subgroup of all $G_i$, in this case $N_1$ is the
intersection of all $N_i$ and any nonzero element in $N_1$ is
primitive.

{\bf Question 1:} Is every irreducible $G$-module quasi-finite
(with respect to a certain quasi-finite sequence of $G$)?

When the irreducible module $N$ is finite dimensional, the answer is
affirmative, since the map $k G=\cup_ikG_i\to\text{End}_kN$ is
surjective and End${}_kN$ is finite dimensional. A weak version of
the above question is the following.

{\bf Question 2:} Assume that $N$ is an irreducible $G$-module. Does
there exist an irreducible $G_i$-submodule $N_i$ of $N$ for each $i$
such that $N$ is the union of all $N_i$.

{\it In the rest of this section $k$ has characteristic 0.}

\subsection{} For quasi-finite groups, a partial
generalization of Mackey's criterion on irreducibility is stated as
follows.

(a) Let $G$ be a quasi-finite group and $H$ a subgroup of $G$. Let
$M$ be a $k H$-module. Then $\ind_H^GM$ is irreducible $G$-module
if  the following two conditions are satisfied.

(1) $M$ is quasi-finitely irreducible (with respect to the
quasi-finite sequence of $H$ induced from the given quasi-finite
sequence of $G$).

(2) Let $M_1$, $M_2$, $M_3,$ ... be a quasi-finite sequence of $M$.
For any positive integer $i$ and $s\in G_i-H\cap G_i$, the two
representations $M_{i,s}$ and $M_i$ of $H_{s,i}=sHs^{-1}\cap H\cap
G_i$ have no common composition factors, where $M_i$ is regarded as
$H_{s,i}$-module by restriction and $M_{i,s}$ is the
$H_{s,i}$-module with $M_i$ as base space and the action of $g\in
H_{s,i}$   on $M_s$ is the same action on $M$ of $s^{-1}gs$.

Proof. Assume the conditions (1) and (2) are satisfied. By Mackey's
criterion, we know that $kG_i\otimes_{kH_i}M_i$  is irreducible
$G_i$-module. By Lemma 1.7 (d) and Lemma 1.6 (b) we see that
$\ind_H^GM$ is irreducible.

\subsection{} Let $A$ be a normal subgroups of a group $G$.
Then for any representation $\rho: A\to GL(V)$ and $s\in G$, we can
define a new representation ${}^s\rho:A\to GL(V)$ by setting
${}^s\rho(g)=\rho(s^{-1}gs)$ for any $g\in A$. In this way we get an
action of $G$ on the set of representations of $A$.

Now assume that that (1) $A$ is commutative and each element of
$A$ has finite order, (2) $G=H\ltimes A$ for some subgroup $H$ of
$G$. By Lemma 1.4, any irreducible representation of $A$ is one
dimensional. Note that the set $X=\Hom(A,k^*)$ is a group. We have
seen that $H$ acts on $X$. Denote by $X/H$ the set of $H$-orbits
in $X$. Let $(\chi_\alpha)_{\alpha\in X/H}$ be a complete set of
representatives of the $H$-orbits. For each $\alpha\in X/H$, let
$H_\alpha$ be the subgroup of $H$ consisting of $h\in H$ with
${}^h\chi_\alpha=\chi_\alpha$ and let $G_\alpha=AH_\alpha$. Define
$\chi_\alpha(gh)=\chi_\alpha(g)$ for any $g\in A$ and $h\in
H_\alpha$. In this way the representation $\chi_\alpha$ is
extended to a representation of $G_\alpha$, denoted again by
$\chi_\alpha$.

Let $\rho$ be an irreducible representation of $H_\alpha$. Through
the homomorphism $G_\alpha\to H_\alpha$ we get an irreducible
representation $\tilde \rho$ of $G_\alpha$. The tensor product
$\tilde\rho\otimes\chi_\alpha$ then is an irreducible
representation of $G_\alpha$. Let
$\theta_{\alpha,\rho}=$Ind${}_{G_\alpha}^G(\tilde\rho\otimes\chi_\alpha)$.

\subsection{Proposition.} Assume that $G$ is quasi-finite. Keep the notations in 1.10.
If  $\rho$ is quasi-finite (with respect to the quasi-finite
sequence of $H_\alpha$ induced from the given quasi-finite
sequence of $G$), then

(a) $\theta_{\alpha,\rho}$ is irreducible.

(b) If $\theta_{\alpha,\rho}$ is isomorphic to
$\theta_{\alpha',\rho'}$, then $\alpha=\alpha'$,  $\rho$ and $\rho'$
are isomorphic.

%(3) If further $H$ is abelian, then any irreducible representation
%of $G$ is isomorphic to a certain $\theta_{\alpha,\rho}$. Moreover
%$\theta_{\alpha,\rho}$ is isomorphic to $\theta_{\alpha',\rho'}$ if
%and only if $\alpha=\alpha'$ and $\rho=\rho'$.

Proof. The argument is similar to that for [S, Proposition 25 (a),
(b)]. Let $G_1,G_2,G_3,...$ be the quasi-finite sequence of $G$.
Then every $G_i$ is finite, $G$  is the union of all $G_i$ and for
any pair $i,j$ there exists $r$ such that both $G_i$ and $G_j$ are
contained in $G_r$. Set $H_{\alpha,i}=H_\alpha\cap G_i$. Let $M$ be
the $kH_\alpha$-module affording the representation $ \rho$ of
$H_\alpha$ and $M_1,M_2,M_3,...$ be a quasi-finite sequence of $M$
(with respect to the sequence
$H_{\alpha,1},H_{\alpha,2},H_{\alpha,3},$...). Let $V$ be the one
dimensional $kG_\alpha$-module affording representation
$\chi_\alpha$. Let $A$ act on each $M_i$ trivially, then $M_i$
becomes an irreducible $G_{\alpha,i}=AH_{\alpha,_i}$-module.
Regarding $V$ as a $G_{\alpha,i}$-module by restriction, then
$M_i\otimes V$ is an irreducible $G_{\alpha,i}$-module.

We claim that Ind${}_{G_{\alpha,i}}^{G_iA}(M_i\otimes V)$ is
irreducible $G_iA$-module. For any $t$ in $G_iA-H_{\alpha,i}A$,
there exists an
 $s$ in $G_i\cap H-G_i\cap H_\alpha$ such that ${}^t\chi_\alpha={}^s\chi_\alpha$. For $s$ in $G_i\cap H-G_i\cap H_\alpha$,
we have ${}^s\chi_\alpha\ne\chi_\alpha$. This implies that there
exists some $a_s\in A$ such that
$\chi_\alpha(a_s)\ne\chi_{\alpha}(s^{-1}a_ss)$. Since $G$ is the
union of all $G_j$ and for any pair $j,j'$ there exists $r$ such
that both $G_j,G_{j'}$ are contained in $G_r$, we can find an $r$
such that $a_s$ is in $G_r$ for any $s$ in $ G_i\cap H-G_i\cap
H_\alpha$, thanks to all $G_j$ being finite. Thus $a_s$ is in $A_{r}=A\cap G_r$ for all $s$ in $
G_i\cap H-G_i\cap H_\alpha$. For $s$ in $ G_i\cap H-G_i\cap
H_\alpha$, set $K_s=H_{\alpha,i}A_r\cap sH_{\alpha,i}A_rs^{-1}$.
Note that the restriction to $H_{\alpha,i}A_r$ of the $H_{\alpha,i}A$-module $M_i\otimes V$ is  irreducible. Through the two injections $K_s\to H_{\alpha,i}A_r$,
$x\to x$ and $x\to s^{-1}xs$ we get two $K_s$-module structures on
the vector space $M_i\otimes V$. The restriction  of the first
$K_s$-module structure  on $M_i\otimes V$ to $A_r$ is the direct sum
of some copies of Res${}_{A_r}\chi_\alpha$, and the restriction of
the second $K_s$-module structure on $M_i\otimes V$ to $A_r$ is the
direct sum of some copies of Res${}_{A_r}{}^s\chi_\alpha$. Since
$a_s$ is in $A_r$ and $\chi_\alpha(a)\ne \chi_\alpha(s^{-1}a_ss)$,
the restrictions  of the two $K_s$-modules to $A_r$ are not
isomorphic, hence the two $K_s$-modules are not isomorphic. By
Mackey's criterion on irreducibility, we see that
Ind${}_{H_{\alpha,i}A_r}^{G_iA_r}(M_i\otimes V)$ is irreducible
$G_iA_r$-module. The natural map
Ind${}_{H_{\alpha,i}A_r}^{G_iA_r}(M_i\otimes
V)\to$Ind${}_{G_{\alpha,i}}^{G_iA}(M_i\otimes V)$ is homomorphism of
$G_iA_r$-module, hence Ind${}_{G_{\alpha,i}}^{G_iA}(M_i\otimes V)$
is irreducible $G_iA$-module. Using Lemma 1.7 (d) and Lemma 1.6 (b) we see that
$\theta_{\alpha,\rho}$ is irreducible.

(2) The restriction of $\theta_{\alpha,\rho}$ to $A$ is completely
reducible and involves only characters in the orbit $H\chi_{\alpha}$
of $\chi_\alpha$, this shows that $\theta_{\alpha,\rho}$ determines
$\alpha$. Let $N$ be the subspace of Ind${}^G_{G_\alpha}(M\otimes
V)$ consisting of all $x\in$Ind${}^G_{G_\alpha}(M\otimes V)$ such
that $\theta_{\alpha,\rho}(a)x=\chi_\alpha(a)x$ for all $a\in A$.
The subspace $N$ is stable under $H_\alpha$ and one checks easily
that the representation of $H_\alpha$ in $N$ is isomorphic to
$\rho$, hence $\theta_{\alpha,\rho}$ determines $\rho$.

The proposition is proved.

\medskip

{\bf Remark.} (1) The above proposition and argument are valid
even if $A$ is not commutative. I am grateful to a referee for
this observation.

(2) It is not clear that whether any irreducible representation of
$G$ is isomorphic to a certain $\theta_{\alpha,\rho}.$

\subsection{} Let  $G$ be a quasi-finite group. Assume that
there exists a sequence $\{1\}=G_0\subset G_1\subset \cdots\subset
G_n=G$ such that all $G_i$ are normal subgroups of $G$ and
$G_i/G_{i-1}$ are abelian. Examples of such groups include Borel subgroups of a reductive group over
$\bbbf$.

\noindent{\bf Question:}  Is each irreducible representation of
$G$  isomorphic to the induced representation of a one dimensional
representation of a subgroup of $G$?

 \section{Algebraic groups with split $BN$-pairs}

\subsection{} In this section we assume that $G$ is an algebraic group with a split $BN$-pair. By definition (see for example [C, p.50]), $G$ has closed subgroups $B$ and $N$ with the following properties,

 (i) The set $B\cup N$ generate $G$, while $T=B\cap N$ is a normal subgroup of $N$ and all elements of $T$ are semisimple .

  (ii) The group $W=N/T$ is generated by a set $S$ of elements $s_i,\ i\in I$, of order 2.

  (iii) If $n_i\in N$ maps to $s_i\in S$ under the natural homomorphism $N\to W$, then $n_iBn_i\ne B$.

  (iv) For each $n\in N$ and each $n_i$ we have $n_iBn\subseteq Bn_inB\cup BnB$.

  (v) $B$ has a closed normal unipotent subgroup $U$ such that $B=T\ltimes U$.

  (vi) $\bigcap_{n\in N}nBn^{-1}=T$

  It is known that $W$ is a  Weyl group. Let $R$ be the root system of $W$ and $\alpha_i,\ i\in I$ are simple roots. For any $w\in W$, $U$ has two subgroups $U_w$ and $U'_w$ such that $U=U'_wU_w$ and $wU'_xw^{-1}\subseteq U$. If $w=s_i$ for some $i$, we simply write $U_i$ and $U'_i$ for $U_w$ and $U'_w$ respectively. For each $w\in W$ we choose an element $n_w\in N$ such that its natural image in $W$ is $w$ and let $n_i$ stand for $n_{s_i}$. The Bruhat decomposition says that $G$ is a disjoint union of the double cosets $Bn_wB,\ w\in W$. Note that  $G_i=B\cup Bn_iB$ is a subgroup $G$.

  Any representation of $T$ can be regarded naturally as a representation of $B$ through the homomorphism $B\to T$. Let $k$ be a field. In this section all representations are assumed over $k$. Let $\theta$ be a one-dimensional representation of $T$, we use the same letter when it is regraded as a representation of $B$. Let $k_\theta$ denote the corresponding $B$-module. We are interested in the induced module $M(\theta)=k G\otimes_{k B}k_\theta.$

  Let $P\supseteq B$ be a parabolic subgroup of $G$ and $L$ a Levi subgroup of $P$ containing $T$. Let $U_P$ be the unipotent radical of $P$. Then $P=L\ltimes U_P$. Moreover $B_L=B\cap L$ is a Borel subgroup and $(B_L,N\cap L) $ forms a $BN$-pair of $L$. By abusing notation, we also use $k_\theta$ for its restriction to $B_L$. Set $M_L(\theta)=kL\otimes_{B_L}k_\theta.$ Let $U_P$ act on $M_L(\theta)$ trivially, then $M_L(\theta)$ becomes a $P$-module. The folowing result is easy to check.

  \subsection{Lemma.} $M(\theta)$ is isomorphic to $kG\otimes_{kP}M_L(\theta)$.

  \medskip

  If $\theta$ is trivial we shall use $M(tr)$ for $M(\theta)$ and $k_{tr}$ for $k_\theta$ respectively.
  Let $1_{tr}$ be a nonzero element in $k_{tr}$. For $x$ in $k G$ we simply denote the element $x\otimes 1_{tr}$ in $M(tr)$ by $x1_{tr}$. For any element $t\in T$ and $n\in N$ we have $nt1_{tr}=n1_{tr}$, so for $w=nT\in W$, the notation $w1_{tr}=n1_{tr}$ is well defined.

  For any subset $J$ of $S$, we shall denote by $W_J$ the subgroup of $W$ generated by $J$ and let $w_J$ be the longest element of $W_J$. Set $\eta_J=\sum_{w\in W_J} (-1)^{l(w)}w1_{tr}$, where $l(w)$ is the length of $w$.
   The following result is a natural extension of [St, Theorem 1, p.348] and part (a) seems new even for finite groups.

 \subsection{Proposition.}   Keep the notations above. Let $J$ be a subset of $S$. Then

 (a) The space $k UW\eta_J$ is a submodule of $M(tr)$ and is denoted by $M(tr)_J$.

 (b)  In particular, $k U\eta_S=kU\sum_{w\in W}(-1)^{l(w)}w1_{tr}$ is a submodule of $M(tr)$. This submodule will be called a Steinberg module of $G$ and is denoted by St.

 Proof. The argument for [St, Theorem 1] works well here.  Clearly $k UW\eta_J$ is stable under the action of $B$. Since $G$ is generated by $B$ and $N$, it remains to check that $k UW\eta_J$ is stable under the action of $N$. But $N$ is generated by all $n_i$ and $T$, so we only need to check that $n_ik UW\eta_J\subseteq k UW\eta_J$. We need to show that $n_iuh\eta_J\in k UW\eta_J$ for any $u\in U$ and $h\in W$. Let $u=u'_iu_i$, where $u_i\in U_i$ and $u'_i\in U'_i$. Then $n_iuh\eta_J=n_iu'_in_i^{-1} n_iu_ih\eta_J.$ Since $n_iu'_in_i^{-1}\in U$, it suffices to check that $n_iu_ih\eta_J\in k UW\eta_J$. When $u_i=1$, this is clear. Now assume that $u_i\ne 1$.  Since $s\eta_J=-\eta_J$ for any $s\in J$, it is no harm to assume that $l(hw_J)=l(h)+l(w_J)$.

 If $hw_J\le s_ihw_J$, then $hw\le s_ihw$ for all $w\in W_J$. In this case, we have $n_iu_ih\eta_J=n_ih\eta_J\in k UW\eta_J$.

 If $s_ih\le h$, then $n_iu_ih\eta_J=n_iu_in_i(s_ih)\eta_J$.
 Note that $n_i^2\in T$. Since $G_i=B\cup Bn_iB$ is a subgroup $G$, if $u_i\ne 1$, we have $n_iu_in_i=n_iu_in_i^{-1}n_i^2=xn_ity$ for some $x,y\in U_i$ and $t\in T$. Thus $n_iu_ih\eta_J=xn_iy(s_ih)\eta_J=xn_i(s_ih)\eta_J=xh\eta_J$ since $(s_ih)w_J\le hw_J$.

 Now assume that $h\le s_ih$ but $s_ihw_J\le hw_J$. Then we must have $s_ih=hs_j$ for some $s_j\in J$. If $w\in W_J$ and $w^{-1}(\alpha_j)$ is a positive root, then we have $hw\le s_ihw$, hence

 (i) $n_iu_ihw1_{tr}=n_ihw1_{tr}=s_ihw1_{tr}=hs_jw1_{tr}$,

 (ii) $hw1_{tr}=xhw1_{tr}$,

 (iii) $n_iu_ihs_jw1_{tr}=n_iu_in_i^{-1}hw1_{tr}=xs_ihw1_{tr}=xhs_jw1_{tr}.$

 Multiplying (i), (ii) and (iii) by $(-1)^{l(w)},\ (-1)^{l(w)}$ and $(-1)^{l(s_jw)}$ respectively, add them, then sum on all $w\in W_J$ satisfying $l(s_jw)=l(w)+1$, then we get
 $$(1-x+n_iu_i)h\eta_J=0.$$
 Thus $n_iuh\eta_J=n_iu'_iu_ih\eta_J=n_iu'_in^{-1}_i n_iu_ih\eta_J=n_iu'_in_i^{-1}(x-1)\eta_J\in kUW\eta_J$. The proposition is proved.

 An analogus of [DL, Prop. 7.3] is the following result.

 \subsection{Proposition.} Let $\theta$ be  a one dimensional representation of $T$, then $M(\theta)\otimes$St is isomorphic to Ind${}_T^Gk_\theta$.

 Proof. Let $1_\theta$ be a nonzero element in $k_\theta$ and $\eta=\sum_{w\in W}(-1)^{l(w)}w1_{tr}$. Then it is easy to check that the map $ g1_\theta\to g(1_\theta\otimes\eta)$ defines an isomorphism of $G$-module between Ind${}_T^Gk_\theta$ and $M(\theta)\otimes$St. The proposition is proved.

 \subsection{Lemma.}  For each $n\in N=N_G(T)$, $kn1_\theta$ is $T$-stable. If each $U_i$ is infinite, then any $T$-stable one dimensional subspace of $M(\theta)$ is contained in $\sum_{n\in N}kn1_\theta$, which is of  dimension $|W|$.

Proof. It is clear.

 \subsection{} Let $J$ be a subset of $S$ and let $M(tr)'_J$ be the sum of all $M(tr)_K$ (see 2.3 (a) for definition) with $J\subsetneq K$. Then $M(tr)'_J$ is a proper submodule of $M(tr)_J$. Let $E_J=M(tr)_J/M(tr)'_J$.

 \subsection{Proposition.} Asusme that each $U_i$ is infinite. If $J$ and $K$ are different subsets of $S$,
 then $E_J$ and $E_K$ are not isomorphic.

 Proof. For any $w\in W$, let $c_w=\sum_{y\le w}(-1)^{l(y)}P_{y,w}(1)y1_{tr}$,
 where $P_{y,w}$ be Kazhdan-Lusztig polynomials. Note that $c_w=\eta_J$ if $w=w_J$ for some subset
 $J$ of $S$.

  We claim that $M(tr)_J$ is the sum of all $kUc_w, \ w\in W$ with $l(ww_J)=l(w)-l(w_J)$. Since $M(tr)_J=kUW\eta_J=kUWc_{w_J}$, we only need to show that $kWc_{w_J}$ is spanned by all $c_w, \ w\in W$ with
  $l(ww_J)=l(w)-l(w_J)$. But this follows from formulas (2.3.a),
  (2.3.c) and Proposition 2.4 in [KL].

   Let $\bar c_w$ be the
 image  of $c_w$ in $E_J$. Let $A_J$ be the subset of $W$ consisting of all $w\in W$ such that $w\le ws$
 for all $s\in S-J$ and $ws\le w$ for all $s\in J$. Then $\bar c_w$ is nonzero if and only if $w\in A_J$ and
 $E_J$ is the  sum of all $kU\bar c_w$. Since $U_i$ is infinite for each $i$, any $T$-stable line in $E_J$ is
 contained in $\sum_{w\in A_J}k\bar c_w=E_J^T$. If there existed an $G$-isomorphism $\phi:E_J\to E_K$, then
 we must have $\phi(E_J^T)=E_K^T$. Thus $\phi(\bar c_{w_J})=\sum_{w\in A_K}a_w\bar c_w,$ $a_w\in k$.
 But $\bar c_{w_J}\ne 0$ is uniquely determined by the following two conditions:
 (1) $n_i\bar c_{w_J}=-\bar c_{w_J}$ if and only if $s_i\in J$,  $U_i\bar c_{w_J}=\bar c_{w_J}$ if and only if
 $s_i\not\in J$. Therefore, $J\ne K$ implies that any nonzero element in $E_K^T$ does not satisfy the conditions for
 $\bar c_{w_J}$, hence $\phi$ does not exist. The proposition is proved.

 \subsection{}  Let $P$ be a parabolic subgroup of $G$ with unipotent radical $U_P$. Assume that $L$ is a Levi subgroup of $P$. Any  $kL$-module $E$ is naturally a $kP$-module through the homomorphism $P\to L$. Then we can define the induced module Ind${}_P^GE=k G\otimes _{kP}E.$ If $P$ contains $B$ and $E$ is one dimensional  $P$-module, then Ind${}_P^GE=k G\otimes _{kP}E$ is a quotient module of some $M(\theta)$.

 Let $P_J$ ($J\subset S$) be a standard parabolic subgroup of $G$. Let $P_J$ act on $k$ trivially. Define $1_{P_J}^G=kG\otimes_{P_J}k$. Clearly $1_{P_J}^G$ is a quotient module of $M(tr)$. Assume that each $U_i$ is infinite. By the discussion above we see that $\Hom_G(1^G_{P_J},1^G_{P_K})$ is nonzero if and only if $J$ is a subset of $K$.

 \section{Reductive Groups with Frobenius Maps}

 \def\bbf{\mathbb{F}_q}
 \def\bbfa{\mathbb{F}_{q^a}}
 \def\bbbf{\bar{\mathbb{F}}_q}
 \def\bbc{\mathbb{C}}

 \subsection{} In this section we assume that $G$ is a connected reductive group defined
 over a finite field $\mathbb{F}_q$ of $q$ elements, where $q$ is a power of a prime $p$.
 Lang's theorem implies that $G$ has a Borel subgroup $B$ defined over $\bbf$ and $B$ contains
 a maximal torus $T$ defined over $\bbf$. For any power $q^a$ of $q$,
 we denote by $G_{q^a}$ the $\mathbb{F}_{q^a}$-points of $G$ and shall identify $G$ with its $\bar{\mathbb{F}}_q$-points,
 where $\bbbf$ is an algebraic closure of $\bbf$. Then we have $G=\bigcup_{a=1}^{\infty}G_{q^a}$.
  Similarly we define $B_{q^a}$ and $T_{q^a}$.

  Let $N$ be the normalizer of $T$ in $G$. Then $B$ and $N$ form a $BN$-pair of $G$.
  Let $k$ be a field. For any one dimensional representation $\theta$ of $T$ over $k$.
  As in section 2 we can define the $k G$-module $M(\theta)=k G\otimes_{k B}k_\theta$. When $\theta$ is trivial representation of $T$ over $k$, as in section 2 we write $M(tr)$ for $M(\theta)$ and let $1_{tr}$ be a nonzero element in $k_\theta$. We shall also write $x1_{tr}$ instead of $x\otimes 1_{tr}$ for $x\in kG$. Let $U$ be the unipotent radical of $B$.

  Recall  that for $w\in W=N/T$, the element $w1_{tr}$ is defined to be  $n_w1_{tr}$, where $n_w$ is a representative in $N$ of $w$  (cf. the paragraph below Lemma 2.2)
.

 \subsection{Theorem.} (a) Assume that $k=\mathbb{C}$ is the field of complex numbers. Then $k U\sum_{w\in W}(-1)^{l(w)}w1_{tr}$ is an irreducible $G$-module.

 (b) Assume that $k=\bbbf$. Then $k U\sum_{w\in W}(-1)^{l(w)}w1_{tr}$ is an irreducible $G$-module.

 Proof. (a) Let $U_{q^a}$ be the $\bbfa$-points of $U$. Then $U=\bigcup_{a=1}^{\infty}U_{q^a}$.
 Let $\eta=\sum_{w\in W}(-1)^{l(w)}w1_{tr}$. Then  $\mathbb{C}[U_{q^a}]\eta$ is isomorphic to the
 Steinberg module of $G_{q^a}$, so it is an irreducible $G_{q^a}$-module. We have
 $\mathbb{C} U\eta=\bigcup_{a=1}^{\infty}\mathbb{C}[U_{q^a}]\eta$ and
 $\mathbb{C}[U_{q^a}]\eta\subset \mathbb{C}[U_{q^{ab}}]\eta$ for any integer $b\ge 1$. By Lemma 1.6 (b),  $\bbc U\eta$ is an irreducible $G$-module.

 The argument for (b) is similar. The theorem is proved.

\subsection{} According to Theorem 2 and Theorem 3 in [St], the $G_{q^a}$-module $k[U_{q^a}]\sum_{w\in W}(-1)^{l(w)}w1_{tr}$ is irreducible if and only if char$k$ does not divide $\sum_{w\in W}q^{al(w)}$. Therefore $k U\sum_{w\in W}(-1)^{l(w)}w1_{tr}$ is irreducible $G$-module if  char$k$ does not divide $\sum_{w\in W}q^{al(w)}$ for all positive integers $a$. Unfortunately it is by no means easy to determine the prime factors of $\sum_{w\in W}q^{al(w)}$ even for type $A_1$ (in this case $W$ has only two elements). So it seems need to find other ways to see whether $k U\sum_{w\in W}(-1)^{l(w)}w1_{tr}$ is
 irreducible $G$-module if  char$k$ is different from 0 and from
 char$\bbbf=p$.

  Let $\theta$ be a group homomorphism from $T$ to $k^*$. For any $w\in W$,
  define ${}^w\theta:T\to k,\ t\to \theta(w^{-1}tw)$.

\subsection{Theorem.}  Assume that $k=\bbc$. Then  $M(\theta)$ has at most $|W_\theta|$ composition factors, where $W_\theta=\{w\in W\,|\, {}^w\theta=\theta\}$. In particular, if  ${}^w\theta\ne \theta$ for any $1\ne w\in W$ (i.e., there exists $t\in T$ such that $\theta (w^{-1}tw)\ne \theta(t)$), then $M(\theta)$ is an irreducible $G$-module.

%(b) $M(\theta)$ and $M(\theta')$ are isomorphic if and only if $\theta'={}^w\theta$ for some $w\in W$.

 Proof.  We can find integer $a$ such that for
any $b\ge a$ we have $W_\theta=\{w\in W\,|\,
{}^w\theta_{T_{q^b}}=\theta_{T_{q^b}}\}$, where $\theta_{T_{q^b}}$
denotes the restriction of $\theta$ to $T_{q^b}$.  Assume that
$0=M_0\subsetneq M_1\subsetneq M_2\subsetneq\cdots\subsetneq
M_h=M(\theta)$ is a filtration of submodules of $M(\theta)$. Then
there exist $x_i\in M_i-M_{i-1}$ for $i=1,2,...,h$. Clearly there
exists $c\ge a$ such that all $x_i$ are in $\bbc G_{q^c}(1\otimes
1_\theta)$, where $1_\theta$ is nonzero element in $k_\theta$. But it is known that $\bbc G_{q^c}(1\otimes 1_\theta)$
has at most $W_\theta$ composition factors.

%(b) Clearly $M(\theta)$ and $M(\theta')$ are isomorphic if  $\theta'={}^w\theta$ for some $w\in W$. The only if part follows from the fact that a $T$-stable line in $M(\theta)$ has the form $n_w\bbc_\theta$, where $n_w$ is a representative in $N$ of $w$.

The theorem is proved.

\subsection{Proposition.} Let $\theta:T\to \bbc^*$ be a group homomorphism. Assume that $W_\theta$ is a parabolic subgroup $W_J$ of $W$. Then the element $\sum_{w\in W_J}(-1)^{l(w)}w1_\theta$ generates an irreducible submodule of $M(\theta)$ and the elements $(s-e)1_{\theta}$, $s\in W_J$  being simple reflections, generate a maximal submodule of $M(\theta)$, where $e$ is the neutral element of $W$.

Proof. It is known that the $kG_{q^a}$-submodule of $M(\theta)$ generated by $\sum_{w\in W_J}(-1)^{l(w)}w1_\theta$ is irreducible for all positive integers $a$ and the $kG_{q^a}$-submodule of $kG1_\theta$ generated by all $(s-e)1_{\theta}$, $s\in W_J$ being simple reflections, is a maxiaml submodule of $kG_{q^a}1_\theta$. The proposition then follows Lemma 1.6 (b).

\subsection{} Assume that $k=\bbc$. It is an interesting question to determine the composition factors of $M(\theta)$.
Assume that $P$ is a parabolic subgroup containing $B$. Let $P$ act trivially on $\bbc$.
Then Ind${}_P^G\bbc=k G\otimes _{kP}\bbc$ is a quotient module of $M(tr)$, so it has finitely many composition factors.
If $P$ is a maximal parabolic subgroup, then Ind${}_P^G\bbc$ has much less composition factors than $M(tr)$.

Let $G$ be a connected reductive group over $\bbf$ such that its
derived group is of type $A_n$. Let $P$ be a maximal parabolic
subgroup of $G$ containing the $F$-stable Borel subgroup $B$ and
assume that the derived subgroup of a Levi subgroup of $P$ has type
$A_{n-1}$. Using 1.6 (b) and representation theory for $G_{q^a}$, it
is easy to see that Ind${}_P^G\bbc$ has a unique irreducible
quotient module which is trivial and a unique  irreducible
submodule.

It is known that there is a bijection between the composition
factors of $G_{q^a}$-submodule $\bbc G_{q^a}1_{tr}$ of $M(tr)$ and
the composition factors of the regular module $\bbc W$ of $W$, which
preserves multiplicities. But this result can not be extended to
$M(tr)$ since by the proof for Proposition 2.7 it is easy to see
that $M(tr)$ has at least $2^{|S|}$ composition factors which are
pairwise non-isomorphic.

\subsection{} Assume that $k=\bbbf$. Then for each dominant
weight $\lambda:T\to k^*$, we have Weyl module $V(\lambda)$ and its
irreducible quotient $L(\lambda)$. Clearly $V(\lambda)$ is a
quotient module of $M(\lambda)$. Also it is clear that the tensor product
$M(\theta)\otimes V(\lambda)$ has a filtration of submodules such
that the quotient modules of the filtration are some $M(\theta+\mu)$, where
$\mu$ are weights of $V(\lambda)$. It is not clear whether some  $M(\theta)$ have infinite composition factors. It might be interesting to study
St$\otimes L(\lambda)$.

It is easy to see that the trivial module of $\bbbf U$ is the
unique irreducible $\bbbf U$-module. A questions comes naturally.
Is every irreducible $\bbbf B$-module one dimensional?

 If  char$k$ is different from 0 and
from
 char$\bbbf=p$, the structure of the modules $M(\theta)$ are more complicated.

\section{Gelfand-Graev Modules}

\subsection{} Keep the notations in 3.1. Thus  $G$ is a
connected reductive group defined over $\bbf$, $B$ a maximal Borel
subgroup of $G$ defined over $\bbf$ and $T$ a maximal torus in $B$
defined over $\bbf$. Let $U$ be the unipotent radical of $B$. The
group $G$ and its subgroups  are identified with their
$\bbbf$-points, so $G=G(\bbbf),\ B=B(\bbbf), \ T=T(\bbbf),$ etc..

Let $R$ be the root system of $G$ and
$\Delta=\{\alpha_1,...,\alpha_l\}$ be the set of simple roots
corresponding to $B$. Denote by $R^+$ the set of positive roots.
For each positive root $\alpha\in R$, let $U_\alpha$ be the
corresponding root subgroup in $G$. We choose an isomorphism
$\varepsilon_\alpha:\bbbf\to U_\alpha$ so that
$t\varepsilon_\alpha(a)t^{-1}=\varepsilon_\alpha(\alpha(t)a)$ for
any $a\in\bbbf$ and $t\in T$. It is known that the subgroup $U'$
of $B$ generated by all $U_\alpha,\ \alpha\in R^+-\Delta$, is a
normal subgroup of $B$ and the quotient group $U/U'$ is isomorphic
to the direct product $U_{\alpha_1}\times
U_{\alpha_2}\times\cdots\times U_{\alpha_l}$ and $B/U'$ is
isomorphic to the semidirect product $T\ltimes U/U'$.

Each irreducible representation of $B/U'$ gives rise naturally an
irreducible representation of $B$. In general it is hard to get a
classification of irreducible representations for groups $B$ and
$U$.

\subsection{} Clearly there is a bijection between
one-dimensional representations of $U$ with $U'$ in the kernel and
the sets $(\sigma_j)$, where $\sigma_j$ is a one-dimensional
representation of $U_{\alpha_j}$. A one-dimensional representation
of $U$ with $U'$ in its kernel is called nondegenerate if all
$\sigma_j$ are non trivial. The group $T$ acts naturally on the set
of irreducible representations of $U/U'$.

It is known that all nondegenerate one-dimensional complex
representations of  $U_q$ form a $T_q$-orbit if the center of $G$
is connected. But the set $\Phi$ of nondegenerate one-dimensional
complex representations of $U$ is uncountable. This implies that
the $T$-orbits in $\Phi$ is uncountable. For a nondegenerate
one-dimensional complex representation $\sigma$ of $U$, we may
consider the induced representation $\ind_U^G\sigma$ and called it
a Gelfand-Graev representation of $G$. It seems not easy to
decompose these Gelfand-Graev representations (cf. section 10 in
[DL])

 \section{Some questions}

 There are some natural questions.

 1. Develop a theory of $k G$-modules for infinite quasi-finite groups. A particular question is to classify irreducible $k G$-modules for some interesting quasi-finite groups, say, reductive groups with Frobenius maps and their Borel subgroups, the infinite Coxeter group $W_{\infty}$ (see 1.8 Example (1)) and the group  $G_{\infty}$ (see 1.8 Example (4)), etc.. .

 According to  Theorem 10.3 and Corollary 10.4 in [BT]) we know that except the trivial representation, all other irreducible representations of $k G$
are infinite dimensional if $G$ is a semisimple algebraic
group over $\bbbf$ and $k$ is infinite with characteristic different
from char$\bbbf$.

 2.  Let $G$ be a connected reductive group over $\bbf$. Then $G$ has a
Frobenius map $F:G\to G$. So for any representation $\rho$ of $G$ we
can define a new representation ${}^F\rho$ by setting
${}^F\rho(g)=\rho(F(g))$. We say that $\rho$ is $F$-stable if $\rho$ and ${}^F\rho$ are isomorphic.

{\bf Question:} Are there any good relations  between the set of
isomorphism classes of irreducible complex representations of $G$
which are $F$-stable and the set of isomorphism classes of
irreducible complex representations of $G^F$?

Replacing $G$ by $GL_n(\bbfa)$ or $SL_n(\bbfa)$,  the above question is answered positively by the theory of Shintani decent (see for example [Sh, Sho, Bo]). For character sheaves, there is a similar result (see [L3]).

 \section{Type $A_1$}

 In this section $G$ will denote  $GL_2(\bbbf)$ or
   $SL_2(\bbbf)$. Let $T$ be the torus of $G$ consisting of diagonal matrices and $B$ the Borel subgroup consisting of
   upper triangle matrices. Let $U$ be the unipotent radical of $B$. In this section we consider complex representations of these groups.

\subsection{} We first consider representations of $B$ over $\bbc$. We have $B=T\ltimes U$. Let $X=\Hom(U,\bbc^*)$. For $t\in T,\ \chi\in X$, define ${}^t\chi:U\to \bbc^*$, $u\to \chi(t^{-1}ut)$. Then we get an action of $T$ on $X$. Note that $U$ is isomorphic to the additive group $\bbbf$, so as abelian group, $U$ is a direct sum of countable
   copies of $\mathbb{F}_p$, where $p$ is the characteristic of $\bbbf$. Therefore, the set $X$ of homomorphism $U\to\bbc^*$ is uncountable.
   This implies the set of $T$-orbits in $X$ is uncountable.

   Denote by $X/T$ the set of $T$-orbits in $X$ and let $(\chi_\alpha)_{\alpha\in X/T}$ be a complete set of representatives of the $T$-orbits. For each
$\alpha\in X/T$, let $T_\alpha$ be the subgroup of $T$ consisting of
$t\in T$ with ${}^t\chi_\alpha=\chi_\alpha$ and let
$B_\alpha=T_\alpha U$. Define $\chi_\alpha(tu)=\chi_\alpha(t)$ for
any $t\in T$ and $u\in U$. In this way the representation
$\chi_\alpha$ is extended to a representation of $B_\alpha$, denoted
again by $\chi_\alpha$. Note that $T_\alpha$ is the center of $B$ if $\chi_\alpha$ is non-trivial, is the whole $T$ if $\chi_\alpha$ is trivial.

Let $\rho$ be an irreducible complex representation of $T_\alpha$, which is one dimensional since $T_\alpha$ is commutative. Through
the homomorphism $B_\alpha\to T_\alpha$ we get an irreducible
representation $\tilde \rho$ of $B_\alpha$. The tensor product
$\tilde\rho\otimes\chi_\alpha$ then is an irreducible representation
of $G_\alpha$. Let $\theta_{\alpha,\rho}$ be the corresponding
induced representation of $B$. According to Proposition 1.11 we have the following result.

\subsection{Lemma.} The complex representation $\theta_{\alpha,\rho}$
of $B$ is irreducible. Moreover
$\theta_{\alpha,\rho}$ is isomorphic to $\theta_{\alpha',\rho'}$ if
and only if $\alpha=\alpha'$ and $\rho=\rho'$.

We can further induce $\theta_{\alpha,\rho}$ to $G$. By the lemma above, if $\chi_\alpha$ is trivial, then $B_\alpha=B$ and $\theta_{\alpha,\rho}$ is just $\tilde\rho$. Since the commutator group $[B,B]$ of $B$ is $U$, any homomorphism $\theta:B\to \bbc^*$ has the form $\tilde\rho$. According to Theorems 3.2 and  3.4, we have the following result.

\subsection{Proposition.} Let $\theta:B\to \bbc^*$ be a group homomorphism. Then

 (a) $M(\theta)=\bbc G\otimes_{\bbc B}\bbc_\theta$ is irreducible $G$-module if $\theta$ is not trivial.

 (b) If $\theta$ is trivial, then $M(\theta)$ has a unique nonzero proper submodule and unique quotient module. The nonzero proper submodule is the Steinberg module. The quotient module is the trivial module of $G$.

 \subsection{} Let ${B_{q}},\ T_{q},\ U_{q}$ be the $F_{q}$-points
   of $B,\ T,\ U$ respectively. Keep the notations in 6.1.  Assume that $\chi_\alpha$ is nontrivial. If the restriction of $\chi_\alpha$ to $U_q$ is not trivial, we can consider the induced representation $\theta_{\alpha,\rho,q}$ of $B_q$ from the restriction of $\tilde\rho\otimes\chi_\alpha$ to $G_\alpha\cap B_q$, which is irreducible. It is known that the action of $B_q$ on $\theta_{\alpha,\rho,q}$ can be extended to actions of $G_q$ and in this way one can get all cuspidal representations of $G_q$. But the author has not been able to see how to extend the $B$ action on $\theta_{\alpha,\rho}$ to an action of $G$.

\section{Miscellany}

In this section we give some discussion to representations of the
groups listed in 1.8 Example (1) and (4).

\subsection{} Let $W=W_\infty$ be the group defined in 1.8 Example (1). Since
$W$ is a Coxeter group, we can use Kazhdan-Lusztig cells to construct
representations of $W$ and its Hecke algebras. Let $s_1,...,s_n$ be
the simple reflections of $W_n$ and let $S$ be the set of all simple
reflections.

(1) Assume that  $W$ is of type $A$.  Let $C_w,\ w\in W$, be the
Kazhdan-Lusztig basis of $\bbc W$ (cf. [KL, Theorem 1.1]). For each
left cell $\sigma$ of $W$, let $I_\sigma$ be the subspace of
$\bbc W$ spanned by all $C_w,\ w\in\sigma$. Denote by $I_{<\sigma}$ the subspace of $\bbc W$ spanned by all $C_w,$ $w\le_L u$ for some $u\in \sigma$ but $w\not\in\sigma$. Then $\bbc W$ is the direct
sum of all $I_\sigma$, and both $I_\sigma+I_{<\sigma}$, $I_{<\sigma}$ are left ideals of $\bbc W$. According to [KL, Theorem 1.4] and 1.6 (b),
 $L_\sigma=(I_\sigma+I_{<\sigma})/I_{<\sigma}$ is an irreducible $\bbc W$-module.  When two left cells $\sigma$ and $\tau$ are in the same two-sided cell, the right star actions leads to an isomorphism between $L_\sigma$ and $L_\tau$. Moreover, $L_\sigma$ and $L_\tau$ are isomorphic
$\bbc W$-modules if and only if  $\sigma$ and $\tau$
are in the same two-sided cell of $W$. Similar results hold for
Hecke algebra of $W$ over $\bbc(q^{\frac12})$ with parameter $q$
(here $q$ is an indeterminate).

According to the proof of [KL, Theorem 1.4], any two-sided cell of
$W$ contains some $w_P$, where  $P$ is a finite subset of $S$ and
$w_P$ is the longest element of the subgroup of $W$ generated by
$P$. Let $\sigma_P$ be the left cell of $W$ containing $w_P$. For subsets of $P$ and $Q$ of $S$, $w_P$ and $w_Q$ are in the
same two-sided cell of $W$ if and
only if $w_P$ and $w_Q$ are in the same two-sided cell of $W_n$
whenever both $w_P$ and $w_Q$ are contained in $W_n$.

Unlike the group $W_n$, some irreducible $\bbc W$-modules are not isomorphic to any  $L_\sigma$, for example, the sign representation of $W$ is not isomorphic to any $L_\sigma$. It is also less easy to discuss irreducible $\bbbf W$-modules.

(2)  Assume that  $W$ is of type $B$.  Let $H$ be the Hecke algebra of $W$ defined over $\mathscr{A}=\bbc[q^{\frac12},q^{-\frac12}]$ ($q$ an indeterminate) with  $\mathscr{A}$-basis $T_w,\ w\in W$, and multiplication relations $(T_{s_i}-q_i)(T_{s_i}+1)=0$ and $T_wT_u=T_{wu}$ if $l(wu)=l(w)+l(u)$, where $q_1=q, \ q_i=q^2$ for all $i\ge 2$. Let $C_w,\ w\in W$, be the
Kazhdan-Lusztig basis of $H$ defined in  [L1,
Proposition 2]). The corresponding cells are called generalized cells
($\varphi$-cells in [L1]. Regarding $\bbc$ as an $\mathscr{A}$-module by specifying $q$ to 1, then we have  $\bbc W=H\otimes_{\mathscr{A}}\bbc$. By abuse notation we use also notation $C_w$ for its image in $\bbc W$.  For each generalized left cell $\sigma$ of
$W$, let $I_\sigma$ be the subspace of
$\bbc W$ spanned by all $C_w,\ w\in\sigma$. Denote by $I_{<\sigma}$ the subspace of $\bbc W$ spanned by all $C_w,$ $w\le_L u$ for some $u\in \sigma$ but $w\not\in\sigma$. Then $\bbc W$ is the direct sum of all $I_\sigma$, and both $I_\sigma+I_{<\sigma}$, $I_{<\sigma}$ are left ideals of $\bbc W$.
According to [L1, Theorem 11] and 1.6 (b), $L_\sigma=(I_\sigma+I_{<\sigma})/I_{<\sigma}$  is an irreducible
$\bbc W$-module.  However, it seems not clear that whether $L_\sigma$
and $L_\tau$ are isomorphic $\bbc W$-modules when the
generalized left cells $\sigma$ and $\tau$ are in the same
generalized two-sided cell of $W$. Similar results hold for the Hecke
algebra $\bar H=H\otimes_{ \mathscr{A}}\bbc(q^{\frac12})$. According to [L1, section 10], if one chooses $q_1=q^3, $ $q_i=q^2$ for all $i\ge 2$, the above results remain valid.

(3) Assume that $W$ is of type $D$. Let $C_w,\ w\in W$, be the
Kazhdan-Lusztig basis of $\bbc W$ (cf. [KL, Theorem 1.1]). For each
left cell $\sigma$ of $W$, as in (1) we can define the subspaces $I_\sigma$ and $I_{<\sigma}$ of $\bbc W$ and the $\bbc W$-module $L_\sigma$. Unlike the case of type $A$ or $B$, $L_\sigma$ may not be irreducible. In [L2,Chapter 12], Lusztig has proved that the $\bbc W_n$-module afforded by a left cell of $W_n$ is multiplicity free and the number of irreducible components in the $\bbc W_n$-module is a power of 2. So it is likely that the $\bbc W$-module $L_\sigma$ is semisimple (that is, a direct sum of some irreducible submodules).

\subsection{} In the rest of this section $G_n$ and  $G_\infty$ are as  in 1.8 Example  (4).
Let $T$ be the subgroup of $G$ consisting of diagonal matrices in $G$ and $N$ be the normalizer of $T$ in $G$.
We can choose naturally a subgroup $B$ of $G$ so that $B$ and $N$ form a $BN$-pair for $G$.
For example, $B$ can be chosen to be the subgroup of $G$ consisting of upper triangular matrices in $G$
if $G=GL_\infty$ or $SL_\infty$. Let $U$ be the kernel of the natural homomorphism $B\to T$. It is esay to
see that $W=N/T$  is just a group in 1.8 Example (1). Let $S=\{s_1,s_2,s_3,...\}$ be the set of simple reflections of $W$.
For each $s_i$ we choose a representative $n_i\in N$ of $s_i$. For each $i$, there exists subgroups $U_i$ and $U_i'$ of
$U$ such that $U=U_i'U_i$ and $n_iU_i'n^{-1}_i\in U$. Nota that if $u_i\in U_i$, we have $n_iu_in_i^{-1}=xn_ity$ for
some $x,y\in U_i$ and $t\in T$.

Set $T_n=T\cap G_n$ and $B_n=B\cap G$. Let $N_n=N\cap G_n$ and
$W_n=N_n/T_n$. Note that $W_n$ can be regarded  as a subgroup of $W$
in a natural way.

Assume that $k=\bar k$,
by Lemma 1.6 (b) the natural representation $V$ of $G$ is
irreducible. We may consider to decompose the tensor product of $m$
copies of $V$. Many classical results for $G_n$ can be extended to
$G_\infty$.

Let $\lambda: T\to k^*$ be a character of $T$. Assume that the
restriction $\lambda_n$  of $\lambda$ to $T_n$ is a dominant weight
of $G_n$ for each $n$. Then we have an irreducible rational
$G_n$-module $V_n$ with highest weight $\lambda_n$. Clearly we have
a natural embedding $V_n\hookrightarrow V_{n+1}$ for each $n$.
Moreover, the embedding is a $G_n$-homomorphism. Let $V(\lambda)$ be
the union of all $V_n$. Then $V(\lambda)$ is naturally a $G$-module.
By Lemma 1.6 (b) it is an irreducible $G$-module.

It is known that a $G_n$-module $M_n$ is called rational if for any
$x\in M$, the $G$-submodule of $G$ generated by $x$ is a finite
dimensional rational $G$-module. We call a $G$-module $M$ rational
if the restriction  of $M$ to $G_n$ is rational. Clearly $V(\lambda)$
 a rational $G$-module in this sense.

\subsection{} Keep the notations in section 7.2. For any homomorphism $\theta:T\to K$, where $K$ is a field, let $K_\theta$ be the corresponding one-dimensional $B$-module. As in section 2 we may consider the induced module $M(\theta)=K G\otimes_{K B}K_\theta$. We can define Steinberg module for $K G$ but which is not a submodule of $M(tr)$.

Let St=$K U\xi$ be a free $K U$-module generated by the  element $\xi$. Note that $K G$ (resp. $k U)$) is the union of all $KG_n$ (resp. $K(U\cap G_n)$). By the proof for Proposition 2.3 we get the following result.

\subsection{Proposition.} The $K U$-module structure on St can be uniquely
extended to a $K G$-module structure as follows:

(1) $tu\xi=tut^{-1}\xi$ for any $t\in T$ and $u\in U$,

(2) $n_iu\xi=-n_iun_i^{-1}\xi$ if $u\in U_i'$,

(3) $n_iu_i'u_i\xi=n_iu_i'n_i^{-1}(x-1)\xi$ for $u_i'\in U_i'$ and $1\ne u_i\in U_i$, where $x\in U_i$ is defined uniquely by the formula $n_iu_in_i^{-1}=xn_ity$, $t\in T,\ y\in U_i$.

Naturally we call the $G$-module St a Steinberg module of $G$. Proposition 2.4 also has its counterpart here, i.e.,  $M(\theta)\otimes$St is isomorphic to Ind${}_T^GK_\theta$.

Using Theorem 2 in [St], Theorem 3.2, Theorem 3.4 and Lemma 1.6 (b) we get the following result.

\subsection{Theorem.} (a) Assume that (1) $k=\bbf$ or $\bbbf$, (2) $K=\bbc$
or $\bbbf$, then St is  irreducible $K G$-module.

(b) Assume that $K=\bbc$ and $\theta:T\to\bbc^*$ is a character of $T$. If $W_\theta=\{w\in W\,|\, {}^w\theta=\theta\}$ has only one element, then $M(\theta)$ is irreducible $K G$-module.

\bigskip

{\bf Acknowledgement.} I thank J. Humphreys and G. Lusztig for
helpful comments. I am grateful to the referees for carefully
reading and helpful comments.

\bigskip

%\newpage

\end{document}